\documentclass{lms}

\usepackage{color}
\usepackage{graphicx,epsf}
\usepackage{amssymb}

\newtheorem{theorem}{Theorem}[section] 
\newtheorem{lemma}[theorem]{Lemma}     
\newtheorem{corollary}[theorem]{Corollary}
\newtheorem{proposition}[theorem]{Proposition}

\newnumbered{assertion}{Assertion}    
\newnumbered{conjecture}{Conjecture}  
\newnumbered{definition}{Definition}
\newnumbered{hypothesis}{Hypothesis}
\newnumbered{remark}{Remark}
\newnumbered{note}{Note}
\newnumbered{observation}{Observation}
\newnumbered{problem}{Problem}
\newnumbered{question}{Question}
\newnumbered{algorithm}{Algorithm}
\newnumbered{example}{Example}
\newunnumbered{notation}{Notation} 

\newcommand{\Real}{{\mathbb R}}

\newcommand{\Z}{\mathbb Z}
\newcommand{\R}{{\rm  R}}

\newcommand{\N}{{\mathbb N}}

\newcommand{\eps}{\varepsilon}

\newcommand{\x}{\mathbf{x}}

\newcommand{\y}{\mathbf{y}}

\newcommand{\z}{\mathbf{z}}

\newcommand {\RR} {{\mathcal R}}

\newcommand {\s}    {{\mbox{\rm sign}}}
\newcommand {\ZZ} {{\rm Z}}

\newcommand {\hide}[1]{}

\title[Homotopy types of fibres]
 {On the number of homotopy types of fibres of a definable map}

\author{Saugata Basu \and Nicolai Vorobjov}

\classno{14P10, 14P25, 32B20}

\extraline{The first author was supported in part by an NSF Career Award 0133597 and
a Alfred P. Sloan Foundation Fellowship.
The second author was supported in part by
the European RTN Network RAAG (contract HPRN-CT-2001-00271).}

\begin{document}
\maketitle

\section{Introduction}
\label{sec:intro}
Let $S \subset \Real^k$ be a set definable in an o-minimal
structure over the reals (see \cite{Dries}), and $\pi_S: S \to
\Real^n$ be a definable map.
For example, $S$ can be a semi-algebraic or a
restricted sub-Pfaffian \cite{GV04} set.
For the purpose of this paper, we will assume without any loss of
generality that $\pi_S$ is the
restriction
to $S$ of the projection map
$\pi:\> \Real^{m+n} \to \Real^n$, where $k=m+n$.

The following statement is a version of Hardt's triviality theorem.

\begin{theorem}[\cite{Hardt,Coste}]
\label{the:hardt}
There exists a finite partition of $\Real^n$ into definable sets
$\{T_i\}_{i \in I}$ such that $S$ is definably trivial over each $T_i$.
\end{theorem}

Theorem \ref{the:hardt} implies
that for each $i \in I$ and any point
$\y \in T_i$, the pre-image $\pi_S^{-1}(T_i)$ is definably homeomorphic to
$\pi_S^{-1}(\y) \times T_i$ by a fiber preserving homeomorphism.
In particular, for each $i \in I$, all fibers $\pi_S^{-1}(\y), \y \in T_i$
are definably homeomorphic.

Hardt's theorem is a corollary of
the existence of {\em cylindrical cell decompositions} of definable sets.
Since the decompositions can be effectively computed in semi-algebraic and
restricted sub-Pfaffian cases (see \cite{BPRbook,GV04}), this implies
a double exponential (in $mn$) upper bound on the cardinality of $I$ and
hence on the number of homeomorphism types of the fibres of the map $\pi_S$.
Apparently, no better bounds than the double exponential bound are known, even though
it seems reasonable to conjecture a single exponential upper bound
on the number of homeomorphism types of the fibres of the map $\pi_S$.

In this paper, we consider the weaker problem of bounding the number of
distinct {\em homotopy types}, occurring amongst the set of all
fibres of $\pi_S$.
The main results of the paper are {\em single exponential} upper bounds on this
number in cases when
$S$ is semi-algebraic, and when $S$ is semi-Pfaffian.
Our results on semi-algebraic sets are formulated and proved with the ground
field being any (possibly non-archimedean) real closed field $\R$.
In the Pfaffian setting we assume that $\R = \Real$.

Even though the precise nature of the bounds we prove is different in
each of the above two cases, the proofs are quite similar.
We will first give a proof in the semi-algebraic case, and then provide
additional details necessary for the Pfaffian case.

The rest of the paper is organized as follows. In Section~\ref{sec:sa}
we state our main result in the semi-algebraic case.
In Section~\ref{sec:proof} the main theorem is proved.
In Section~\ref{sec:pfaffian} we state the result in the Pfaffian
case, and discuss the modifications needed in its proof.
In Section~\ref{sec:fewnomials} we use the result of Section~\ref{sec:pfaffian} to state
and prove a single exponential upper bound on the number of homotopy
types of semi-algebraic sets defined by fewnomials,
as well as those defined by polynomials with bounded additive complexity.
Finally, in
Section~\ref{sec:metric} we prove some metric upper bounds, related to
homotopy types, for semi-algebraic sets defined by polynomials with
integer coefficients.

\section{Semi-algebraic case}\label{sec:sa}

Let $\R$ be a real closed field,
${\mathcal P} \subset \R[X_1,\ldots,X_m,Y_1,\ldots,Y_n]$, and
let $\phi$ be a Boolean formula with atoms of the form
$P=0$, $P > 0$, or $P< 0$, where $P \in {\mathcal P}$.
We call $\phi$ a ${\mathcal P}$-formula, and the semi-algebraic set
$S \subset \R^{m+n}$ defined by $\phi$, a ${\mathcal P}$-semi-algebraic set.
Note that for a given ${\mathcal P}$ there is a finite number of different
${\mathcal P}$-semi-algebraic sets.
Also, it is clear from the definition that
every semi-algebraic set is a ${\mathcal P}$-semi-algebraic set
for an appropriate ${\mathcal P}$.

If the Boolean formula $\phi$ contains no negations, and
its atoms are of the form
$P= 0$, $P \geq 0$, or $P \leq  0$, with $P \in {\mathcal P}$,
then we call $\phi$ a ${\mathcal P}$-closed formula, and the semi-algebraic set
$S \subset \R^{m+n}$ defined by ${\phi}$,
a ${\mathcal P}$-closed semi-algebraic set.

\subsection{Main Result}\label{subsec:mainresult}

We prove the following theorem.

\begin{theorem}
\label{the:main}
Let ${\mathcal P} \subset \R[X_1,\ldots,X_m,Y_1,\ldots,Y_n]$,
with $\deg(P) \leq d$ for each $P \in {\mathcal P}$ and cardinality
$\#{\mathcal P} = s$.
Then, there exists a finite set $A \subset \R^n$,
with \[
\# A \leq (2^m snd)^{O(nm)},
\]
such that for every $\y \in \R^n$ there exists $\z \in A$ such that
for every ${\mathcal P}$-semi-algebraic set $S \subset \R^{m+n}$, the set
$\pi_S^{-1} (\y)$
is semi-algebraically homotopy equivalent to
$\pi_S^{-1} (\z)$.
In particular, for any fixed ${\mathcal P}$-semi-algebraic set $S$,
the number of different homotopy types of fibres $\pi_{S}^{-1} (\y)$
for various $\y \in \pi(S)$ is also bounded by
\[
(2^m snd)^{O(nm)}.
\]
\end{theorem}

\begin{remark}
We actually prove a more precise bound,
$$\# A \le s^{2(m+1)n} (2^m nd)^{O(nm)}$$
(see Section~\ref{subsec:proof}).
\end{remark}

Notice that the bound in Theorem~\ref{the:main} is single exponential in $mn$.
The following example shows that the single exponential dependence on $m$ is unavoidable.

\begin{example}
Let
$P \in \R[X_1,\ldots,X_m] \hookrightarrow \R[X_1,\ldots,X_m,Y]$
be the polynomial defined by
\[
P := \sum_{i=1}^{m} \prod_{j=1}^{d}(X_i - j)^2.
\]
The algebraic set defined by $P=0$ in
$\R^{m+1}$
with coordinates $X_1,\ldots,X_m,Y$,
consists of $d^m$ lines all parallel to the $Y$ axis.
Consider now the semi-algebraic set $S \subset \R^{m+1}$ defined by
$$
\displaylines{
(P = 0)\; \wedge \;
(0 \le Y \le X_1+ dX_2 + d^2 X_3 + \cdots + d^{m-1}X_m).
}
$$
It is easy to verify that, if $\pi:\> \R^{m+1} \to \R$ is the
projection map on the $Y$ co-ordinate, then
the fibres $\pi_S^{-1}(y)$, for $y \in \{0,1, 2, \ldots ,d^m-1 \} \subset \R$
are 0-dimensional and of different cardinality, and hence have different homotopy types.
\end{example}

The above example does not exhibit exponential dependence of the
number of homotopy types of fibres on $n$,
which is equal to $1$ in the example.
In fact, we cannot hope to produce examples where the number of homotopy types
of the fibres grows with $n$ (with the parameters $s,d,$ and $m$ fixed)
since this number can be bounded by a function of $s,d$ and $m$
independent of $n$, as we show next.

Suppose that $\phi$ is a Boolean formula with atoms
$\{a_i,b_i,c_i \mid 1 \leq i \leq s\}$.
For an ordered list ${\mathcal P} = (P_1,\ldots,P_s)$ of polynomials
$P_i \in {\R[X_1,\ldots,X_m]}$, we denote by $\phi_{\mathcal P}$
the formula obtained from $\phi$ by replacing
for each $i,\ 1\leq i \leq s $, the atom $a_i$
(respectively, $b_i$ and $c_i$) by
$P_i= 0$ (respectively, by $P_i > 0$ and by $P_i < 0$).

\begin{definition}\label{def:hom_lists}
We say that two ordered lists ${\mathcal P} = (P_1,\ldots,P_s)$,
${\mathcal Q} = (Q_1,\ldots,Q_s)$
of polynomials $P_i, Q_i \in {\R[X_1,\ldots,X_m]}$ have the same {\em homotopy
type} if for any Boolean formula $\phi$, the semi-algebraic sets defined by
$\phi_{\mathcal P}$ and $\phi_{\mathcal Q}$ are homotopy equivalent.
\end{definition}

Let
${\mathcal S}_{m,s,d}$ be the family of all ordered lists $(P_1,\ldots,P_s)$,
$P_i \in \R[X_1,\ldots,X_m]$, with $\deg(P_i) \leq d$ for $1 \leq i \leq s$.

\begin{corollary}
\label{cor:main}
The number of different homotopy types of ordered lists in ${\mathcal S}_{m,s,d}$
does not exceed
\begin{equation}\label{form:main}
\left( s {m+d \choose d} \right)^{O \left( s{m+d \choose d}m \right)}=(sd)^{O(sd^m)}.
\end{equation}
In particular, the number of different homotopy types of semi-algebraic sets
defined by a fixed formula $\phi_{\mathcal P}$, where ${\mathcal P}$ varies over
${\mathcal S}_{m,s,d}$, does not exceed (\ref{form:main}).
\end{corollary}

\begin{proof}
We introduce
one new variable for each
of the possible ${m+d \choose d}$ monomials of degree at most $d$ in
$m$ variables, so that every polynomial $P_i$
is uniquely represented by a point in $\R^{{m+d \choose d}}$.
Thus, every ordered list in ${\mathcal S}_{m,s,d}$
is uniquely represented by a point in $\R^{s{m+d \choose d}}$.
Consider the ordered list ${\mathcal Q}=(Q_1, \ldots ,Q_s)$ of polynomials
$$Q_i \in \Z[Y_1, \ldots ,Y_{s{m+d \choose d}},X_1, \ldots ,X_m],$$
$1 \le i \le s$, of degrees at most $d+1$, such that for each $\y \in \R^{s{m+d \choose d}}$
the list
$$(Q_1 (\y, X_1, \ldots , X_m), \ldots , Q_s (\y, X_1, \ldots , X_m)) \in
{\mathcal S}_{m,s,d}.$$
Now apply Theorem~\ref{the:main} to ${\mathcal Q}$ and the projection map
$\pi':\> \R^{s{m+d \choose d}+m} \to \R^{s{m+d \choose d}}$.
\end{proof}

Corollary~\ref{cor:main} implies that
for every ${\mathcal P}$-semi-algebraic set $S \subset \R^{m+n}$, where
${\mathcal P} \in {\mathcal S}_{m+n,s,d}$, the number of different
homotopy types of fibres $\pi_{S}^{-1} (\y)$ does not exceed
(\ref{form:main}) (as well as the bound from Theorem~\ref{the:main}).
In particular, this number has an upper bound not depending on $n$.

\section{Proof of Theorem~\ref{the:main}}\label{sec:main}
\label{sec:proof}
\subsection{Main ideas}\label{ideas}
We first summarize the main ideas behind the proof of Theorem~\ref{the:main}.

Let $\R$ be a real closed field and
${\mathcal P} = \{P_1,\ldots,P_s\}
\subset {\R[X_1,\ldots,X_m,Y_1,\ldots,Y_n]}$, with
$\deg(P_i) \leq d, 1 \leq i \leq s$.
We fix a finite set of points
$B \subset \R^n$ such that for every $\y \in \R^n$ there exists $\z \in B$ such that
for every ${\mathcal P}$-semi-algebraic set $S$, the set $\pi_S^{-1}(\y) $
is semi-algebraically homotopy equivalent to $\pi_S^{-1}(\z) $.
The existence of a set $B$ with this property follows from
Hardt's triviality theorem
(Theorem \ref{the:hardt}),
and the fact that the number of ${\mathcal P}$-semi-algebraic sets is finite.

Observe that it is possible to choose
$\omega \in \R, \omega > 0$,
sufficiently large (depending on ${\mathcal P}$ and $B$) such that,
for any ${\mathcal P}$-semi-algebraic set $S$,
\begin{enumerate}
\item
the intersection of $S$  with the set defined by the conjunction
of the $2(m+n)$ inequalities
$-\omega < X_i < \omega$, $-\omega < Y_j < \omega$,
$1 \le i \le m,\ 1 \le j \le n$,
has the same homotopy type as $S$ (denote this intersection by $S_0$);
\item
the set $B$ lies in $(-\omega, \omega)^n$, and $B$ preserves its defining property
with respect to $S_0$, i.e.,
for every $\y \in \R^n$ there exists $\z \in B$ such that the set
$\pi_{S_0}^{-1}(\y) $
is semi-algebraically homotopy equivalent to
$\pi_{S_0}^{-1}(\z) $.
\end{enumerate}
We will henceforth assume that
${\mathcal P}$ contains the $2(m+n)$ polynomials,
$X_i \pm \omega, Y_j \pm \omega , 1\leq i \leq m, 1 \leq j \leq n$, and
restrict our attention to the bounded ${\mathcal P}$-semi-algebraic sets, which are
contained in the cube $(-\omega,\omega)^{m+n}$.
Note that in order to bound the number of homotopy types of the fibres
$\pi_S^{-1}(\y), \y \in \pi(S)$, where $S$ is an arbitrary
${\mathcal P}$-semi-algebraic set,
it is sufficient to bound the number of
homotopy types of the fibres
$\pi_{S_0}^{-1}(\y), \y \in \pi(S_0),$
of the projection of the bounded set $S_0$.

We replace ${\mathcal P}$ by another family
${\mathcal P}' \subset {\R[X_1,\ldots,X_m,Y_1,\ldots,Y_n]}$
having the following properties.

\begin{enumerate}
\item
The cardinality $\#{\mathcal P}' = 2s^2$, and
$\deg(P) \leq d$, for each $P \in {\mathcal P}'$.
\item
For each bounded ${\mathcal P}$-semi-algebraic set $S$, there exists a
bounded and ${\mathcal P}'$-closed semi-algebraic set $S'$, such that
$S' \simeq S$ (where $\simeq$ stands for semi-algebraic homotopy equivalence).
\item
$\pi_{S'}^{-1}(\y) \simeq \pi_S^{-1}(\y)$ for every $\y \in B$.
\item
The zero sets of the polynomials in ${\mathcal P}'$ are non-singular
hypersurfaces, intersecting transversally, and
the same is true for the restriction of ${\mathcal P}'$ on
each  $\pi^{-1}(\y)$, $\y \in B$.
It follows that the partition of $\R^{m+n}$ into ${\mathcal P}'$-sign
invariant subsets is a Whitney stratification, denoted by ${\mathcal W}$, of $\R^{m+n}$,
which is compatible with every ${\mathcal P}'$-semi-algebraic set.
\end{enumerate}

Let $G_1$ be the set of all critical values of $\pi$,
restricted to various strata of the stratification ${\mathcal W}$
of dimensions not less than $n$, and let $G_2$
be the union of the projections of all strata of ${\mathcal W}$ of
dimensions less than $n$.
Let
$G:= G_1 \cup G_2$.
We prove that $B \subset \R^n \setminus G$,
and hence in order to bound the number of homotopy types of fibres $\pi_S^{-1}(\y)$
for any bounded ${\mathcal P}$-semi-algebraic set $S$,
it is sufficient to bound the number of
homotopy types of fibres, $\pi_{S'}^{-1}(\y)$, $\y \in \R^n \setminus G$.
We then prove, using Thom's first isotopy lemma,
that the homeomorphism type of the fibre $\pi_{S'}^{-1}(\y)$ does not change
as $\y$ varies over any fixed connected component of $\R^n  \setminus G$.
This, along with property (3) of ${\mathcal P}'$ stated above, shows that in order to
bound the number of
homotopy types of the fibres $\pi_S^{-1}(\y)$ it suffices to bound the
number of connected components of $\R^n \setminus G$.
We obtain the set $A$ (in Theorem \ref{the:main}) by choosing
one point in each connected component of $\R^n \setminus G.$
In order to bound the cardinality of $A$, we use a recent result from \cite{GVZ}
bounding the Betti numbers of projections of semi-algebraic sets
in terms of those of certain fibred products.

\subsection{Notations}\label{sub:notations}

Let $\R$ be  a real closed field. For an element $a \in \R$ introduce
$$
\s(a) = \Biggl\{
\begin{tabular}{ccc}
0 & \mbox{ if }  $a=0$,\\
1 & \mbox{ if } $a> 0$,\\
$-1$& \mbox{ if } $a< 0$.
\end{tabular}
$$
If ${\mathcal P} \subset {\R} [X_1, \ldots , X_k]$ is finite, we write the
{\em set of zeros} of ${\mathcal P}$ in ${\R}^k$ as
$$
\ZZ({\mathcal P}):=\Bigl\{ \x \in {\R}^k\mid\bigwedge_{P\in{\mathcal P}}P(\x)= 0 \Bigr\}.
$$
A  {\em sign condition} $\sigma$ on
${\mathcal P}$ is an element of $\{0,1,- 1\}^{\mathcal P}$.
The {\em realization of the sign condition
$\sigma$} is the basic semi-algebraic set
$$
\RR(\sigma) := \Bigl\{ \x \in {\R}^k\;\mid\;
\bigwedge_{P\in{\mathcal P}} \s({P}(\x))=\sigma(P) \Bigr\}.
$$
A sign condition $\sigma$ is {\em realizable} if $\RR(\sigma) \neq \emptyset$.
We denote by ${\rm Sign}({\mathcal P})$ the set of realizable sign conditions
on ${\mathcal P}.$
For $\sigma \in {\rm Sign}({\mathcal P})$ we define the {\em level of} $\sigma$
as the cardinality $\#\{P \in {\mathcal P}| \sigma(P) = 0 \}$.
For each level $\ell$, $0 \leq \ell \leq \# {\mathcal P}$, we denote by
${\rm Sign}_{\ell}({\mathcal P})$
the subset of ${\rm Sign}({\mathcal P})$ of elements of level $\ell$.

Finally, for a sign condition $\sigma$ let
$$
{\mathcal Z}(\sigma) := \Bigl\{ \x \in {\R}^k\;\mid\;
\bigwedge_{P\in{\mathcal P},\ \sigma (P)=0} P(\x)=0 \Bigr\}.
$$

\subsection{Replacing a bounded set by a homotopy equivalent closed bounded set}

In this section, we describe a modification
of the construction from \cite{GV} (see also \cite{BPR05})
for replacing any given bounded semi-algebraic set by a
closed bounded  semi-algebraic set
which has the same homotopy type as the original set.

\begin{definition}
Let ${\mathcal F}(x)$ be a predicate defined over $\R_+$ and $y \in \R_+$.
The notation $\forall(0 < x \ll y)\ {\mathcal F}(x)$
stands for the statement
$$\exists z \in (0,y)\ \forall x \in \R_+\ ({\rm if}\> x<z,\> {\rm then}\> {\mathcal F}(x)),$$
and can be read ``for all positive
$x$ sufficiently smaller than $y$, ${\mathcal F(x)}$ is true''.
\end{definition}

\begin{definition}\label{def:small}
For points
$$\bar \eps =(\varepsilon_{2s,1}, \ldots ,
\eps_{2s,s}, \varepsilon_{2s-1,1}, \ldots , \eps_{2s-1,s}, \cdots ,\eps_{1,1}, \ldots ,
\eps_{1,s}) \in \R_{+}^{2s^2},$$
and a predicate ${\mathcal F}(\bar \eps)$ over $\R_{+}^{2s^2}$
we say
``for all sufficiently small $\bar \eps$, ${\mathcal F}(\bar \eps)$ is true'' if
$$
\forall(0 < \eps_{2s,1} \ll 1)
\forall(0 < \eps_{2s,2} \ll \eps_{2s,1})
\cdots \forall(0 < \eps_{1,s} \ll \eps_{1,s-1})
{\mathcal F}(\bar \eps).
$$
\end{definition}

Consider a finite set of polynomials
$${\mathcal P} = \{P_1,\ldots,P_s\} \subset
\R[X_1,\ldots,X_m,Y_1,\ldots,Y_n].$$
and a sequence of elements in $\R$
$$
1 > \varepsilon_{2s,1} > \cdots
> \eps_{2s,s} > \varepsilon_{2s-1,1} > \cdots > \eps_{2s-1,s} > \cdots > \eps_{1,s} >0.
$$
Given $\sigma \in {\rm Sign}_{\ell}({\mathcal P})$ we denote by
$\RR(\sigma_+)$  the closed semi-algebraic set defined by
the conjunction of the inequalities:
\begin{eqnarray*}
&~&- \eps_{2 \ell,i} \le P_{i} \le  \eps_{2 \ell,i},  \mbox{ for each }
P_i \in {\mathcal P}
\mbox{ such that } \sigma(P_i) = 0,\\
&~&P \ge  0,  \mbox{ for each }  P \in {\mathcal P}
\mbox{ such that } \sigma(P) = 1,\\
&~&P \le  0,  \mbox{ for each }  P \in {\mathcal P}
\mbox{ such that } \sigma(P) = -1.
\end{eqnarray*}
We denote by
$\RR({\sigma_-})$
the open semi-algebraic set defined by the conjunction of the inequalities:
\begin{eqnarray*}
&~&- \eps_{2 \ell -1,i} < P_{i} <  \varepsilon_{2 \ell -1,i}, \mbox{ for each }
P_i \in {\mathcal P}
\mbox{ such that } \sigma(P_i) = 0,\\
&~&P >  0,  \mbox{ for each }  P \in {\mathcal P}
\mbox{ such that } \sigma(P) = 1,\\
&~&P <  0,  \mbox{ for each }  P \in {\mathcal P}
\mbox{ such that } \sigma(P) = -1.
\end{eqnarray*}
Notice that,
$\RR(\sigma) \subset \RR(\sigma_+)$ and $\RR(\sigma) \subset \RR(\sigma_-)$.

Now suppose that the bounded ${\mathcal P}$-semi-algebraic set $S$ is defined
by
$$
\displaylines{
S = \bigcup_{\sigma \in \Sigma_S}\RR(\sigma),}
$$
for some $\Sigma_S\subset {\rm Sign}({\mathcal P})$.

\begin{definition}
\label{def:S'}
Let
$$
\displaylines{
\Sigma_{S,\ell}= \Sigma_S \cap {\rm Sign}_{\ell}({\mathcal P})
}
$$
and define a sequence of sets,
$S_{\ell}(\bar\eps) \subset {\R}^{m+n}$, $0 \leq \ell \leq s$ inductively as follows.

\begin{itemize}
\item[$\bullet$]
$S_{0}(\bar\eps) := S$.
\item[$\bullet$]
For
$0 \leq \ell \leq s$,
let
$$
\displaylines{
S_{\ell +1}(\bar\eps) =
\left(
S_{\ell}(\bar\eps) \setminus
\bigcup_{\sigma \in {\rm Sign}_{\ell}({\mathcal P})\setminus \Sigma_{S,\ell}}\RR({\sigma_-})
\right)
\cup \bigcup_{\sigma \in \Sigma_{S,\ell}}\RR({\sigma_+}).
}
$$
\end{itemize}

We denote $S'(\bar\eps) := S_{s+1}(\bar\eps)$ and
${\mathcal P}' := \{ P_{j} \pm \eps_{i,j} \mid 1 \leq i \leq 2s, 1 \leq j \leq s\}$.
\end{definition}

The following statement is proved in \cite{BPR05} and is
a strengthening  of a result in \cite{GV} (where it is shown that
for all sufficiently small $\bar\eps$ the
sum of the Betti numbers of $S$ and $S'(\bar\eps)$ are equal).

\begin{proposition}[\cite{BPR05}]
\label{the:GV}
For all sufficiently small $\bar\eps$ (see Definition~\ref{def:small}),
$S'(\bar\eps) \simeq S$.
\end{proposition}

We now show that it is possible to
rewrite the Boolean formula for $S'(\bar\eps)$, which originally
(in its negation-free form) contains inequalities
of the kind $P \ge 0$ or $P \le 0$ for $P \in {\mathcal P}$,
so that it still defines
$S'(\bar\eps)$
but involves only inequalities of the kind
$Q \ge 0$ or $Q \le 0$, for $Q \in {\mathcal P}'$.

\begin{lemma}
\label{lem:closed}
For all sufficiently small $\bar\eps$,
$S'(\bar\eps)$ is a ${\mathcal P}'$-closed semi-algebraic set.
\end{lemma}

\begin{proof}
In the negation-free Boolean formula for $S'(\bar\eps)$ replace every inequality
of the kind $P_j \ge 0$
(respectively, of the kind $P_j \le 0$), where
$P_j \in {\mathcal P}$, by $P_j \ge \eps_{2,j}$
(respectively, by $P_j \le - \eps_{2,j}$).
Denote by $S''(\bar\eps)$ the set defined by the new formula.
Obviously $S''(\bar\eps)$ is a ${\mathcal P}'$-closed semi-algebraic set and
moreover it is clear from the definition that
$S''(\bar\eps) \subset S'(\bar\eps)$.

We now prove that $S'(\bar\eps) \subset S''(\bar\eps)$.
Let $\z \in S'(\bar\eps)$.
Let $\RR({\sigma_+})$, where $\sigma \in \Sigma_{\ell}$, be the maximal
(with respect to $\ell$) non-empty set in the definition of $S'(\bar\eps)$
containing $\z$.
Assume that
$$\RR({\sigma_+})= \{ (\x, \y)|\>  - \eps_{2 \ell,i_1} \le P_{i_1} \le  \eps_{2 \ell,i_1},
\ldots, - \eps_{2 \ell,i_{\ell}} \le P_{i_{\ell}} \le  \eps_{2 \ell,i_{\ell}},$$
$$P_{i_{\ell +1}} \ge 0, \ldots , P_{i_k} \ge 0 \}.$$
We show that
$$\z \in \{ (\x, \y)|\>  - \eps_{2 \ell,i_1} \le P_{i_1} \le  \eps_{2 \ell,i_1},
\ldots, - \eps_{2 \ell,i_{\ell}} \le P_{i_{\ell}} \le  \eps_{2 \ell,i_{\ell}},$$
$$P_{i_{\ell +1}} \ge \eps_{2,i_{\ell +1}} , \ldots , P_{i_k} \ge \eps_{2, i_k} \},$$
and hence $\z \in S''(\bar\eps)$.
Indeed, suppose that $P_{i_r}(\z) < \eps_{2, i_r}$
for some $r$, ${\ell+1} \le r \le k$.
Then $- \eps_{2 (\ell+1),i_r} \le P_{i_r}(\z) \le  \eps_{2 (\ell+1),i_r}$, since
$0 \le P_{i_r}(\z) < \eps_{2, i_r} \le \eps_{2 (\ell+1),i_r}$.
Therefore $\z \in \RR({\sigma'_+})$, where $\sigma' \in \Sigma_{\ell+1}$,
which contradicts the maximality of $\sigma$ with respect to $\ell$.
\end{proof}
In particular, we have that
$S'(\bar\eps)$ is a ${\mathcal P}'$-semi-algebraic set,
and we define $\Sigma_S' \subset {\rm Sign}({\mathcal P}')$ by
\[
S'(\bar\eps) =\bigcup_{ \sigma \in \Sigma_{S}'} \RR(\sigma).
\]

We now note some extra properties of the family ${\mathcal P}'$.

\begin{lemma}\label{prop:A}
For all sufficiently small $\bar\eps$, the following holds.
If $\sigma \in {\rm Sign}_{\ell}({\mathcal P}')$, then $\ell \le m+n$ and
$\RR(\sigma) \subset {\R}^{m+n}$ is a non-singular $(m+n-\ell)$-dimensional manifold
such that at every point $(\x,\y) \in \RR(\sigma)$, the
$(\ell \times (m+n))$-Jacobi matrix,
\[
\left( \frac{\partial P}{\partial X_i} , \frac{\partial P}{\partial Y_j}
\right)_{P \in{\mathcal P}',\ \sigma(P) = 0,\
1\leq i \leq m,\ 1 \leq j \leq n}
\]
has the maximal rank $\ell$.
\end{lemma}

\begin{proof}
Suppose without loss of generality that
$$\{ P \in {\mathcal P}' |\> \sigma (P)=0 \}= \{ P_{i_1}-\eps_{j_1,i_1},
\ldots , P_{i_{\ell}}-\eps_{j_{\ell},i_{\ell}} \}$$
(cf. definition of ${\mathcal P}'$),
where $1 > \eps_{j_1,i_1} > \cdots > \eps_{j_{\ell},i_{\ell}} >0$.
Let $\ell \le m+n$.
Consider the semi-algebraic
map $P_{i_1,\ldots,i_{\ell}}: {\R}^{m+n} \rightarrow {\R}^{\ell}$ defined by
$$(\x,\y) \mapsto (P_{i_1}(\x,\y),\ldots,P_{i_{\ell}}(\x,\y)).$$
By the semi-algebraic version of Sard's theorem (see \cite{BCR}), the
set of critical values of $P_{i_1,\ldots,i_{\ell}}$ is a semi-algebraic
subset of ${\R}^{\ell}$ of dimension strictly less than $\ell$.
It follows that {\em for all sufficiently small} $\bar\eps$, the point
$(\eps_{j_1,i_1},\ldots,\eps_{j_{\ell},i_{\ell}})$ is not a critical value
of the map $P_{i_1,\ldots,i_{\ell}}$.

It follows that the algebraic set
$$\{ (\x, \y)|\> P_{i_1}=\eps_{j_1,i_1}, \ldots , P_{i_{\ell}}=\eps_{j_{\ell},i_{\ell}} \}$$
is a smooth $(m+n- \ell)$-dimensional manifold
such that at every point on it
the $(\ell \times (m+n))$-Jacobi matrix,
\[
\left( \frac{\partial P}{\partial X_i} , \frac{\partial P}{\partial Y_j}
\right)_{P \in \{ P_{i_1}, \ldots , P_{i_{\ell}}\},\
1\leq i \leq m,\ 1 \leq j \leq n}
\]
has the maximal rank $\ell$.
The same is true for the basic semi-algebraic set $\RR(\sigma)$.

We now prove that $\ell \le m+n$.
Suppose that $\ell > m+n$.
As we have just proved,
$\{ P_{i_1}=\eps_{j_1,i_1}, \ldots , P_{i_{m+n}}=\eps_{j_{m+n},i_{m+n}} \}$
is a finite set of points.
Polynomial $P_{i_{\ell}}-\eps_{j_{\ell},i_{\ell}}$ vanishes on each of these points.
Choosing a value of $\eps_{j_{\ell},i_{\ell}}$ smaller than
the value of $P_{i_{\ell}}$ at any point, we get a contradiction.
\end{proof}

Recall that $B \subset \R^n$ is a fixed finite set of points
such that for every $\y \in \R^n$ there exists $\z \in B$ such that
for every ${\mathcal P}$-semi-algebraic set $S$, we have
$\pi_S^{-1}(\y) \simeq \pi_S^{-1}(\z) $.

\begin{lemma}\label{prop:B}
For every $\y \in B$, bounded ${\mathcal P}$-semi-algebraic set
$S$, for all sufficiently small $\bar\eps$,
and $\sigma \in {\rm Sign}_{\ell}({\mathcal P}_{\y}')$,
where
${\mathcal P}_{\y}':= \{ P(X_1, \ldots ,X_m, \y)|\> P \in {\mathcal P}' \}$,
the following holds.
\begin{enumerate}
\item
$0 \leq \ell \leq m$, and $\RR(\sigma) \cap \pi^{-1}(\y)$
is a non-singular $(m-\ell)$-dimensional manifold
such that at every point $(\x,\y) \in \RR(\sigma) \cap \pi^{-1}(\y)$,
the $(\ell \times m)$-Jacobi matrix,
\[
\left( \frac{\partial P}{\partial X_i} \right)_{P \in{\mathcal P}_{\y}',\ \sigma(P) = 0,\
1 \leq i \leq m}
\]
has the maximal rank $\ell$;
\item
$\pi^{-1}(\y) \cap S'(\bar\eps) \simeq \pi_S^{-1}(\y)$.
\end{enumerate}
\end{lemma}

\begin{proof}
Note that $P_{\y} \in \R[X_1,\ldots,X_m]$ for each $P \in {\mathcal P}$
and $\y \in {\R}^n$.
The proof of (1) is now identical to the proof of Lemma~\ref{prop:A}.
The part (2) of the lemma is a consequence of Proposition~\ref{the:GV}.
\end{proof}

\subsection{Whitney stratification of $\R^n$ compatible with ${\mathcal P}'$}
\label{sec:whitney}

\begin{lemma}\label{Whitney}
For any bounded ${\mathcal P}$-semi-algebraic set $S$ and
for all sufficiently small $\bar\eps$, the partitions
\begin{eqnarray*}
\label{partition}
\R^n &=& \bigcup_{ \sigma \in {\rm Sign}({\mathcal P}')} \RR(\sigma),\\
S'(\bar\eps) &=& \bigcup_{ \sigma \in \Sigma_S'} \RR(\sigma),
\end{eqnarray*}
are  compatible Whitney stratifications of $\R^n$ and $S'(\bar\eps)$ respectively.
\end{lemma}

\begin{proof}
Follows directly from the definition of Whitney stratification (see \cite{GM,CS}),
and Lemma~\ref{prop:A}.
\end{proof}

Fix some sign condition $\sigma \in {\rm Sign}({\mathcal P}')$.
Recall that $(\x,\y) \in \RR (\sigma)$ is a {\em critical point}
of the map $\pi_{\RR (\sigma)}$ if the Jacobi matrix,
\[
\left( \frac{\partial P}{\partial X_i} \right)_{P \in{\mathcal P}',\ \sigma(P) = 0,\
1 \leq i \leq m}
\]
at $(\x, \y)$ is not of the maximal possible rank.
The projection $\pi (\x, \y)$ of a critical point is a {\em critical value}
of $\pi_{\RR (\sigma)}$.

Let $C_1(\bar\eps) \subset \R^{m+n}$
be the set of critical points of $\pi_{\RR (\sigma)}$
over all sign conditions
$$\sigma \in \bigcup_{\ell \le m} {\rm Sign}_{\ell}({\mathcal P}'),
$$
(i.e., over all $\sigma \in {\rm Sign}_{\ell}({\mathcal P}')$
with $\dim (\RR (\sigma)) \ge n$).
For a bounded ${\mathcal P}$-semi-algebraic set $S$,
let $C_1(S,\bar\eps) \subset S'(\bar\eps)$
be the set of critical points of $\pi_{\RR (\sigma)}$
over all sign conditions
$$\sigma \in \bigcup_{\ell \le m} {\rm Sign}_{\ell}({\mathcal P}')\cap
\Sigma_S'$$
(i.e., over all $\sigma \in \Sigma_S'$ with $\dim (\RR (\sigma)) \ge n$).

Let $C_2(\bar\eps) \subset \R^{m+n}$
be the union of $\RR (\sigma)$ over all
$$\sigma \in \bigcup_{\ell > m} {\rm Sign}_{\ell}({\mathcal P}')$$
(i.e., over all $\sigma \in {\rm Sign}_{\ell}({\mathcal P}')$
with $\dim (\RR (\sigma)) < n$).
For a bounded ${\mathcal P}$-semi-algebraic set $S$,
let $C_2(S,\bar\eps) \subset S'(\bar\eps)$
be the union of $\RR (\sigma)$ over all
$$\sigma \in \bigcup_{\ell > m} {\rm Sign}_{\ell}({\mathcal P}') \cap \Sigma_S' $$
(i.e., over all $\sigma \in \Sigma_S'$ with $\dim (\RR (\sigma)) < n$).
Denote
$C(\bar\eps):= C_1(\bar\eps) \cup C_2(\bar\eps)$, and
$C(S,\bar\eps):= C_1(S,\bar\eps) \cup C_2(S,\bar\eps)$.
Notice that,
$C(S,\bar\eps) = C(\bar\eps) \cap S'(\bar\eps).$
\begin{lemma}\label{closed}
For each bounded ${\mathcal P}$-semi-algebraic $S$ and for all sufficiently small $\bar\eps$,
the set $C(S,\bar\eps)$ is closed and bounded in ${\R}^{m+n}$.
\end{lemma}

\begin{proof}
The set $C(S,\bar\eps)$ is bounded since $S'(\bar\eps)$ is bounded.
The union $C_2(S,\bar\eps)$ of strata of dimensions less than $n$ is closed
since $S'(\bar\eps)$ is closed.

Let $\sigma_1 \in {\rm Sign}_{\ell_1}({\mathcal P}') \cap \Sigma_S'$,
$\sigma_2 \in {\rm Sign}_{\ell_2}({\mathcal P}') \cap \Sigma_S'$,
where $\ell_1 \le m$, $\ell_1 < \ell_2$,
and if $\sigma_1 (P)=0$, then $\sigma_2 (P)=0$ for any $P \in {\mathcal P}'$.
It follows that stratum $\RR (\sigma_2)$ lies in the closure of the stratum $\RR (\sigma_1)$.
Let ${\mathcal J}$ be the finite family of $(\ell_1 \times \ell_1)$-minors such that
$Z( {\mathcal J}) \cap \RR (\sigma_1)$ is the set of all critical points of
$\pi_{\RR (\sigma_1)}$.
Then $Z( {\mathcal J}) \cap \RR (\sigma_2)$ is either contained in
$C_2(S,\bar\eps)$
(when $\dim (\RR (\sigma_2)) <n$), or is contained in the set of all critical points
of $\pi_{\RR (\sigma_2)}$ (when $\dim (\RR (\sigma_2)) \ge n$).
It follows that the closure of $Z( {\mathcal J}) \cap \RR (\sigma_1)$ lies in the union
of the following sets:
\begin{enumerate}
\item
$Z( {\mathcal J}) \cap \RR (\sigma_1)$,
\item
sets of critical points of some strata of dimensions less than $m+n- \ell_1$,
\item
some strata of the dimension less than $n$.
\end{enumerate}
Using induction on descending dimensions in case (2), we conclude that the closure of
$Z( {\mathcal J}) \cap \RR (\sigma_1)$ is contained in $C(S, \bar\eps)$.
Hence, $C(S, \bar\eps)$ is closed.
\end{proof}

\begin{definition}
\label{def:criticalvalues}
We denote by
$G_i(\bar\eps) := \pi(C_i(\bar\eps)), i= 1,2$ and
$G(\bar\eps) := G_1(\bar\eps) \cup G_2(\bar\eps)$.
Similarly, for each bounded ${\mathcal P}$-semi-algebraic set $S$,
we denote by
$G_i(S,\bar\eps) := \pi(C_i(S,\bar\eps)), i= 1,2$ and
$G(S,\bar\eps) := G_1(S,\bar\eps) \cup G_2(S,\bar\eps)$.
\end{definition}

\begin{lemma}\label{representatives}
For all sufficiently small $\bar\eps$, $B \cap G(\bar\eps) = \emptyset$.
In particular, $B \cap G(S,\bar\eps) = \emptyset$
for every bounded ${\mathcal P}$-semi-algebraic set $S$.
\end{lemma}

\begin{proof}
By Lemma~\ref{prop:B},
for all $\y \in B$, for all sufficiently small $\bar\eps$,
and $\sigma \in {\rm Sign}_{\ell}({\mathcal P}_{\y}')$,
\begin{enumerate}
\item
$0 \leq \ell \leq m$, and
\item
$\RR(\sigma) \cap \pi^{-1}(\y)$
is a non-singular $(m-\ell)$-dimensional manifold
such that at every point $(\x,\y) \in \RR(\sigma) \cap \pi^{-1}(\y)$,
the $(\ell \times m)$-Jacobi matrix,
\[
\left( \frac{\partial P}{\partial X_i} \right)_{P \in{\mathcal P}_{\y}',\ \sigma(P) = 0,\
1 \leq i \leq m}
\]
has the maximal rank $\ell$.
\end{enumerate}
If a point $\y \in B \cap G_1(\bar\eps) = B \cap \pi(C_1(\bar\eps))$, then
there exists $\x \in \R^m$ such that $(\x,\y)$
is a critical point of $\pi_{\RR (\sigma)}$
for some $\sigma \in \bigcup_{\ell \le m} {\rm Sign}_{\ell}({\mathcal P}')$,
and this is impossible by (2).

Similarly, $\y \in B \cap G_2(\bar\eps) = B \cap \pi(C_2(\bar\eps))$,
implies that there exists $\x \in \R^m$ such that
$(\x,\y) \in \RR (\sigma)$ for some
$\sigma \in \bigcup_{\ell > m} {\rm Sign}_{\ell}({\mathcal P}')$, and this
is impossible  by (1).
\end{proof}

For all sufficiently small $\bar\eps$, let $D(\bar \eps)$ be a connected component of
$\R^n \setminus G(\bar\eps)$, and for a bounded
${\mathcal P}$-semi-algebraic set $S$,
let $D(S,\bar\eps)$ be a connected component of $\pi(S'(\bar\eps)) \setminus G(S,\bar\eps)$.

\begin{lemma}\label{prop:discriminant}
For every bounded ${\mathcal P}$-semi-algebraic set $S$
and for all sufficiently small $\bar\eps$,
all fibers $\pi^{-1}(\y) \cap S'(\bar\eps)$, $\y \in D(\bar\eps)$
are homeomorphic.
\end{lemma}

\begin{proof}
Lemma~\ref{Whitney} implies that
$\widehat S:=\pi_{S'(\bar\eps)}^{-1}(\pi(S'(\bar\eps)) \setminus
G(S,\bar\eps))$ is a
Whitney stratified set having strata of dimensions at least $n$.
Moreover, $\pi_{\widehat S}$ is a proper stratified submersion.
By Thom's first isotopy lemma (in the semi-algebraic version, over real closed fields
\cite{CS}) the map $\pi_{\widehat S}$ is a locally trivial fibration.
In particular, all fibers $\pi_{S'(\bar\eps)}^{-1}(\y)$, $\y \in D(S,\bar\eps)$
are homeomorphic for every connected component $D(S,\bar\eps)$.
The lemma follows, since
the inclusion $G(S,\bar\eps) \subset G(\bar\eps)$ implies that either
$D(\bar \eps) \subset D(S,\bar\eps)$ for some connected component $D(S,\bar\eps)$, or
$D(\bar \eps) \cap \pi(S'(\bar\eps))= \emptyset$.
\end{proof}

\begin{lemma}\label{components}
For each $\y \in B$
and for all sufficiently small $\bar\eps$,
there exists a connected component
$D(\bar\eps)$ of $\R^n \setminus G(\bar\eps)$, such that
$\pi_{S}^{-1}(\y) \simeq \pi_{S'(\bar\eps)}^{-1}(\y_1)$
for every bounded ${\mathcal P}$-semi-algebraic set $S$ and
for every $\y_1 \in D(\bar\eps)$.
\end{lemma}

\begin{proof}
By Lemma~\ref{representatives},
$\y$ belongs to some connected component $D(\bar\eps)$ for all sufficiently small $\bar\eps$.
Lemma~\ref{prop:B}~(2) implies that $\pi^{-1}(\y) \cap S \simeq \pi^{-1}(\y) \cap S'(\bar\eps)$.
Finally, by Lemma~\ref{prop:discriminant},
$\pi_{S'(\bar\eps)}^{-1}(\y) \simeq \pi_{S'(\bar\eps)}^{-1}(\y_1)$
for every $\y_1 \in D(\bar\eps)$.
\end{proof}

\subsection{Proof of Theorem~\ref{the:main}}\label{subsec:proof}

The following proposition gives a bound on the Betti numbers of the
projection $\pi(V)$ of a closed and bounded semi-algebraic set $V$
in terms of the number and degrees of polynomials defining $V$.

\begin{proposition}[\cite{GVZ}]\label{prop:GVZ}
Let $V \subset {\R}^{m+n}$ be a closed and bounded semi-algebraic set defined by a
Boolean formula with $s$ distinct polynomials of degrees not exceeding $d$.
Then the $k$-th Betti number of the projection
$${\rm b}_k (\pi (V)) \le (ksd)^{O(n+km)}.$$
\end{proposition}
\begin{proof}
See \cite{GVZ}.
\end{proof}

For the proof of Theorem~\ref{the:main} we need
the following inequalities which can be derived from the Mayer-Vietoris
exact sequence (see, for instance, Proposition 7.33 in \cite{BPRbook}).
\begin{proposition}[(Mayer--Vietoris inequalities)]\label{prop:MV}
Let the subsets $W_1, \ldots , W_t \subset \R^n$ be all open or all closed.
Then
\begin{equation}\label{eq:MV1}
{\rm b}_k \left( \bigcup_{1 \le j \le t} W_j \right) \le \sum_{J \subset \{1, \ldots ,\ t \}}
{\rm b}_{k- (\# J) +1} \left( \bigcap_{j \in J} W_j \right)
\end{equation}
and
\begin{equation}\label{eq:MV2}
{\rm b}_k \left( \bigcap_{1 \le j \le t} W_j \right) \le \sum_{J \subset \{1, \ldots ,\ t \}}
{\rm b}_{k+ (\# J) -1} \left( \bigcup_{j \in J} W_j \right),
\end{equation}
where ${\rm b}_k$ is the $k$-th Betti number.
\end{proposition}

\begin{proof}[of Theorem~\ref{the:main}]
To make clearer the idea of the proof, we first show a weaker bound,
$\# A \le (s^m nd)^{O(nm)}$, which requires fewer technical details.

Let
$$
{\mathcal P}' = \{Q_1,\ldots,Q_{2s^2}\},
$$
(see Definition~\ref{def:S'} for the definition of the set ${\mathcal P}'$).
Recall that for all sufficiently small $\bar\eps$,
$G(\bar\eps)=G_1(\bar\eps) \cup G_2(\bar\eps)$,
where $G_1(\bar\eps)$ is the union of sets of critical values of
$\pi_{\RR (\sigma)}$ over all strata $\RR (\sigma)$ of dimensions at least $n$,
and $G_2(\bar\eps)$ is the union of projections of all strata of dimensions less than $n$.

Observe that the set of critical points $C_1(\bar\eps)$ is a
$\widehat {\mathcal P}$-semi-algebraic set, where
\[
\widehat {\mathcal P} = {\mathcal P}' \cup {\mathcal Q},
\]
\[
{\mathcal Q} = \left\{ \det \left( M^{i_1,\ldots,i_\ell}_{j_1,\ldots,j_\ell} \right)
\;\mid\; 1 \leq \ell \leq m, 1 \leq i_1 < \cdots < i_\ell \leq 2s^2,
1 \leq j_1 < \cdots < j_\ell \leq m \right\},
\]
and
$$
M^{i_1,\ldots,i_\ell}_{j_1,\ldots,j_\ell}=
\left(
\begin{array}{ccc}
\frac{\partial Q_{i_1}}{\partial X_{j_1}} &\cdots &
\frac{\partial Q_{i_1}}{\partial X_{j_\ell}}\\
\vdots &\vdots & \vdots\\
\frac{\partial Q_{i_\ell}}{\partial X_{j_1}} &\cdots &
\frac{\partial Q_{i_\ell}}{\partial X_{j_\ell}}\\
\end{array}
\right).
$$
Since $\#{\mathcal P}' \le 2s^2$, we conclude that
$\# \widehat {\mathcal P} \le s^{O(m)}$.
The degrees of polynomials from $\widehat {\mathcal P}$ do not exceed $O(md)$.
The union $C_2(\bar\eps)$
of all strata of dimensions less than $n$ is a
${\mathcal P}'$-semi-algebraic set.
It follows that $C(\bar\eps)$ is a $\widehat {\mathcal P}$-semi-algebraic set
which is closed and bounded in ${\R}^{m+n}$ due to Lemma~\ref{closed}.

Then Proposition~\ref{prop:GVZ} implies that
${\rm b}_{n-1}(G(\bar\eps)) \le (s^mnd)^{O(nm)}$.
By Alexander's duality, the homology groups
$\widetilde H_0({\R}^n \setminus G(\bar\eps))$ and $H_{n-1}(G(\bar\eps))$
are isomorphic, and hence the number of connected components
$${\rm b}_0({\R}^n \setminus G(\bar\eps)) \le (s^mnd)^{O(nm)}.$$

We choose for the set $A$, one point in each connected component of
$\R^n \setminus G(\bar\eps)$.
Now the bound $\# A \le (s^m nd)^{O(nm)}$ follows from Lemma~\ref{components}.

Now we proceed to the proof of the bound
$${\rm b}_0({\R}^n \setminus G(\bar\eps)) \le s^{2(m+1)n} (2^m nd)^{O(nm)}$$
which implies the same
bound for $\# A$.

Recall (Section~\ref{ideas}) that the set $S$ is contained in $(- \omega, \omega)^{n+m}$ for
a sufficiently large $\omega >0$, and therefore
$S'(\bar \eps)$ lies in the closed box $[- \omega -1, \omega +1]^{n+m}$.
Denote by ${\mathcal E}_1(\bar \eps)$ the family of closed sets
each of which is the intersection of $[- \omega -1, \omega +1]^{n+m}$ with the set
of critical points of
$\pi_{{\mathcal Z} (\sigma)}$, over all sign conditions $\sigma$ such that
strata $\RR (\sigma)$ have dimensions at least $n$ (the notation ${\mathcal Z} (\sigma)$
was introduced in Section~\ref{sub:notations}).
Let ${\mathcal E}_2(\bar \eps)$ be the family of closed sets
each of which is the intersection of $[- \omega -1, \omega +1]^{n+m}$ with
${\mathcal Z} (\sigma)$, over all
sign conditions $\sigma$ such that strata $\RR (\sigma)$ have dimensions equal to $n-1$.
Let ${\mathcal E}(\bar \eps):= {\mathcal E}_1(\bar \eps) \cup {\mathcal E}_2(\bar \eps)$.
Denote by $E(\bar \eps)$ the image under the projection $\pi$ of the union of all sets in
the family ${\mathcal E}(\bar \eps)$.

Because of the transversality condition,
every stratum of the stratification of $S'(\bar \eps)$, having the dimension
less than $n+m$, lies in the closure of a stratum, having the next higher dimension.
In particular, this is true for strata of dimensions less than $n-1$.
It follows that $G(\bar \eps) \subset E(\bar \eps)$, and thus
every connected component of the complement $\R^n \setminus E(\bar \eps)$
is contained in a connected component of $\R^n \setminus G(\bar \eps)$.
Since $\dim (E(\bar \eps))<n$, every connected component of $\R^n \setminus G(\bar \eps)$
contains a connected component of $\R^n \setminus E(\bar \eps)$.
Therefore, it is sufficient to estimate from above the
Betti number ${\rm b}_0 (\R^n \setminus E(\bar \eps))$ which is equal to
${\rm b}_{n-1}(E(\bar \eps))$ by the Alexander's duality.

The total number of sets ${\mathcal Z} (\sigma)$, such that
$\sigma \in {\rm Sign}({\mathcal P}')$ and $\dim ({\mathcal Z} (\sigma)) \ge n-1$,
is $O(s^{2(m+1)})$ because each ${\mathcal Z} (\sigma)$ is
defined by a conjunction of at most $m+1$ of possible $2s^2$ polynomial equations.
Thus, the cardinality $\# {\mathcal E}(\bar \eps)$, as well as the
number of images under the projection $\pi$ of sets in
${\mathcal E}(\bar \eps)$ is
$O(s^{2(m+1)})$.
According to (\ref{eq:MV1}) in Proposition~\ref{prop:MV},
${\rm b}_{n-1}(E(\bar \eps))$
does not exceed the sum of certain Betti numbers of sets of the type
$$
\Phi :=\bigcap_{1 \le i \le p} \pi (U_i),
$$
where every $U_i \in {\mathcal E}(\bar \eps)$ and $1 \leq p \leq n$.
More precisely, we have
$$
\displaylines{
{\rm b}_{n-1}(E(\bar \eps)) \;\leq \;\sum_{1 \le p \le n}\quad
\sum_{ \{ U_{1}, \ldots ,U_{p} \} \subset\ {\mathcal E}(\bar \eps)}
{\rm b}_{n-p} \left( \bigcap_{1 \le i \le p} \pi (U_i) \right).
}
$$
Obviously, there are $O(s^{2(m+1)n})$ sets of the kind $\Phi$.

Using inequality (\ref{eq:MV2}) in
Proposition \ref{prop:MV}, we have that for each $\Phi$ as above,
the Betti number ${\rm b}_{n-p}(\Phi)$ does not exceed
the sum of certain Betti numbers of unions of the kind,
$$\Psi := \bigcup_{1 \le j \le q} \pi (U_{i_j}) =
\pi \left( \bigcup_{1 \le j \le q} U_{i_j} \right),$$
with  $1 \leq q \leq p$.
More precisely,
\begin{eqnarray*}
{\rm b}_{n-p} (\Phi) &\;\leq\;&
\sum_{1 \le q \le p}\quad \sum_{1 \leq i_1 < \cdots< i_q \leq p}
{\rm b}_{n-p+q-1} \left( \pi \left( \bigcup_{1 \le j \le q} U_{i_j} \right) \right).
\end{eqnarray*}
It is clear that there are at most $2^{p} \leq 2^n$ sets of the kind $\Psi$.

If a set $U \in {\mathcal E}_1(\bar \eps)$, then it is defined by $O(n+m)$ polynomials
of degrees at most $d$ (including the
linear polynomials defining $[- \omega -1, \omega +1]^{n+m}$).
If a set $U \in {\mathcal E}_2(\bar \eps)$, then it is defined by $O(n+2^m)$ polynomials
of degrees $O(md)$,
since the critical points on strata of dimensions at least $n$
are defined by $O(2^m)$ determinantal equations, the
corresponding matrices have orders $O(m)$, and
the entries of these matrices are polynomials of degrees at most $d$.

It follows that the closed and bounded set
$$\bigcup_{1 \le j \le q} U_{i_j}$$
is defined by $O(n(n+2^m))$ polynomials of degrees $O(md)$.
By Proposition~\ref{prop:GVZ},
${\rm b}_{n-p+q-1}(\Psi) \le (2^mnd)^{O(nm)}$ for all $1 \le p \le n$, $1 \le q \le p$.
Then ${\rm b}_{n-p} (\Phi) \le (2^mnd)^{O(nm)}$ for every $1 \le p \le n$.
Since there are $O(s^{2(m+1)n})$ sets of the kind $\Phi$, we get the
claimed bound
$${\rm b}_{n-1}(E(\bar \eps)) \le s^{2(m+1)n}(2^mnd)^{O(nm)}.$$
\end{proof}

\section{Semi-Pfaffian case}\label{sec:pfaffian}

Pfaffian functions, introduced in \cite{Kh}, are real analytic
functions satisfying triangular systems of Pfaffian (first order partial differential)
equations with polynomial coefficients.

\begin{definition}[\cite{Kh,GV04}]
A {\em Pfaffian chain} of the order $r \ge 0$ and degree $\alpha \ge 1$ in
an open domain $U \subset \Real^n$ is a sequence of real analytic functions
$f_1, \ldots , f_r$ in $U$ satisfying differential equations
\begin{equation}\label{pfaff}
df_j(\x)= \sum_{1 \le i \le n} g_{ij}(\x, f_1(\x), \ldots ,f_j(\x))dx_i
\end{equation}
for $1 \le j \le r$.
Here $g_{ij}(\x, y_1, \ldots,y_j)$ are polynomials in $\x=(x_1, \ldots, x_n)$,
$y_1, \ldots, y_j$ of degrees not exceeding $\alpha$.
A function $f(\x)=P(\x,f_1(\x), \ldots ,f_r(\x))$, where the polynomial
$P(\x, y_1,\ldots , y_r)$ has a degree not exceeding
$\beta \ge 1$, is called a {\em Pfaffian function} of order $r$ and
degree $(\alpha, \beta)$.
\end{definition}

The class of Pfaffian functions includes polynomials, real algebraic functions,
elementary transcendental functions defined in appropriate domains (see \cite{GV04}).
One of the important subclasses is formed by {\em fewnomials} \cite{Kh}.

Let ${\mathcal P}$ be a finite set of Pfaffian functions in the open cube
$U:=(-1,1)^{m+n} \subset \Real^{m+n}$.

\begin{definition}
A set $S \subset U$ is called ${\mathcal P}$-semi-Pfaffian in
$U$ (or just {\em semi-Pfaffian}
when some ${\mathcal P}$ is fixed) if it is defined by a
Boolean formula with atoms of the form
$P > 0,\ P < 0,\ P=0$ for $P \in {\mathcal P}$.
A ${\mathcal P}$-semi-Pfaffian set $S$ is {\em restricted} if its closure in
$U$ is compact.
\end{definition}

Semi-Pfaffian sets share many of the finiteness properties of semi-algebraic
sets \cite{GV04}.
The proofs of Proposition \ref{the:GV}  and all lemmas from
Section~\ref{sec:main} extend to the semi-Pfaffian case without any difficulty.

\begin{theorem}\label{the:main2}
Let ${\mathcal P}$ be a finite set of Pfaffian functions defined
on the open cube
$U:=(-1,1)^{m+n} \subset \Real^{m+n}$, with
$\#{\mathcal P} =s$, and such that
all functions in ${\mathcal P}$ have degrees $(\alpha, \beta)$ and
are derived from a common Pfaffian chain of order $r$.
Then, there exists a finite set $A \subset \pi(U)$
with
\[
\#A \leq s^{O(nm)}2^{O(n(m^2+nr^2))}(nm(\alpha+ \beta))^{O(n(m+r))},
\]
such that for every $\y \in \pi(U)$ there exists $\z \in A$ such that
for every ${\mathcal P}$-semi-Pfaffian set $S \subset U$, the set
$\pi_S^{-1} (\y) $
is homotopy equivalent to
$\pi_S^{-1} (\z) $.
In particular, for any fixed ${\mathcal P}$-semi-Pfaffian  set $S$,
the number of different homotopy types of fibres $\pi_{S}^{-1} (\y)$
for various $\y \in \pi(S)$ is also bounded by
\[
s^{O(nm)}2^{O(n(m^2+nr^2))}(nm(\alpha+ \beta))^{O(n(m+r))}.
\]
\end{theorem}

\begin{proof}
We add to ${\mathcal P}$ the $2(m+n)$ functions,
$X_i \pm (1 - \delta), Y_j \pm (1-\delta), 1 \le i \le m,\ 1 \le j \le n$,
where $\delta \in \Real$ is chosen sufficiently small and positive.
We restrict attention to those ${\mathcal P}$-semi-Pfaffian sets
which are contained in the cube defined by,
$-1+ \delta < X_i < 1- \delta$, $-1 + \delta < Y_j < 1 - \delta$,
$1 \le i \le m,\ 1 \le j \le n$.
For any ${\mathcal P}$-semi-Pfaffian set $S$,
the intersection $S_0:= S \cap (-1 + \delta, 1- \delta)^{n+m}$,
is homotopy equivalent to $S$.

As in the proof of Theorem \ref{the:main}, let
$B \subset \Real^n$ is a fixed finite set of points
such that for each  ${\mathcal P}$-semi-Pfaffian set
$S$ and $\y \in \pi(S)$, there exists $\z \in B \cap \pi(S)$ (depending on
$S$) such that
$\pi_{S}^{-1}(\y) \simeq \pi_{S}^{-1}(\z)$.

Then, for all sufficiently small $\delta > 0$,
for each ${\mathcal P}$-semi-pfaffian set $S$,
$\pi_{S_0}^{-1}(\y) \simeq \pi_{S}^{-1}(\y)$ for each $\y \in B$.
It follows that it is sufficient to bound from above the number
of homotopy types
of fibres of the projection of the restricted sets $S_0$.

The rest of the proof of Theorem~\ref{the:main2}
is identical to the corresponding part of
the proof of Theorem~\ref{the:main}.
The only essential difference is the replacement of the reference to
Proposition~\ref{prop:GVZ} by the reference to the following
proposition.

\begin{proposition}[\cite{GVZ}]
Let $V \subset U$ be a closed and bounded ${\mathcal P}$-semi-Pfaffian set defined by a
Boolean formula such that $\#{\mathcal P} =s$,
all functions from ${\mathcal P}$ are defined in $U$, have degrees $(\alpha, \beta)$, and
a common Pfaffian chain is of order $r$.
Then the $k$-th Betti number of the projection
$${\rm b}_k (\pi (V)) \le (ks)^{O(n+km)}2^{O(kr)^2}
((n+km)(\alpha + \beta))^{O(n+km+kr)}.$$
\end{proposition}
\end{proof}

In Section~\ref{sec:fewnomials} we will need the following
technical improvement of Theorem~\ref{the:main2}.
We first introduce some notations.

\begin{definition}\label{def:octant}
Let $\tau \in {\rm Sign}(\{X_1, \ldots ,X_m \})$ be a sign condition
for the family of coordinate functions.
We denote the realization ${\mathcal R}(\tau)$ of $\tau$ by $\Real^{m}_{\tau}$ and call it
an {\em octant} of $\Real^m$.
In general, for any subset $S \subset \Real^m$ and
$\tau \in {\rm Sign}(\{X_1, \ldots ,X_m \})$ we introduce
$S_\tau := S \cap \Real^{m}_{\tau}$.
\end{definition}

\begin{theorem}\label{the:main3}
Let
$${\mathcal P} = \bigcup_{\tau \in {\rm Sign}(\{ X_1, \ldots , X_m \})}
{\mathcal P}_{\tau}$$
be a finite set of Pfaffian functions such that for each
$\tau \in {\rm Sign}(\{ X_1, \ldots , X_m \})$
functions in ${\mathcal P}_{\tau}$
are defined in the domain
$U_{\tau} \times V$, where $U$ is either $(-1,1)^m$ or $\Real^{m}$, and $V$
is either $(-1,1)^n$ or $\Real^{n}$.
Let $\#{\mathcal P} =s$,
all functions in ${\mathcal P}$ have degrees $(\alpha, \beta)$, and
for each $\tau$ functions in ${\mathcal P}_{\tau}$
are derived from a common Pfaffian chain of order at most $r$.
Then, there exists a finite set $A \subset U$ with
$$\# A \le  s^{O(nm)}2^{O(n(m^2+nr^2))}(nm(\alpha+ \beta))^{O(n(m+r))},$$
such that for every $\y \in U$ there exists $\z \in A$ such that
for every set $S= \bigcup_{\tau}S_{\tau} \subset (U \times V)$,
where every $S_{\tau}$ is a ${\mathcal P}_{\tau}$-semi-Pffafian set,
the set $\pi_{S}^{-1}(\y)$ is homotopy equivalent to $\pi_{S}^{-1}(\z)$.
In particular, for any fixed set $S$,
the number of different homotopy types of fibres $\pi_{S}^{-1} (\y)$
for various $\y \in \pi(S)$ is bounded by
\[
s^{O(nm)}2^{O(n(m^2+nr^2))}(nm(\alpha+ \beta))^{O(n(m+r))}.
\]
\end{theorem}

\begin{proof}
The proof is a straightforward modification of the proof of Theorem
\ref{the:main2}.
We express $S$ as a union of realizations of sign conditions $\sigma \wedge \tau$
on the families ${\mathcal P}_\tau \cup \{X_1,\ldots,X_m\}$,
where $\sigma \in {\rm Sign}({\mathcal P}_{\tau})$ and
$\tau \in {\rm Sign}(\{ X_1, \ldots , X_m \})$.
Thus,
\begin{equation}\label{eq:s}
S = \bigcup_{\tau \in {\rm Sign}(\{ X_1, \ldots , X_m \}),\>
 \sigma \in \Sigma_{\tau,S}}\RR(\sigma \wedge \tau)
\end{equation}
with $\Sigma_{\tau, S}\subset {\rm Sign}({\mathcal P}_\tau)$.

Given $\tau \in {\rm Sign}_{\ell}(\{X_1,\ldots,X_m\})$,
let $\tau_>$ (respectively, $\tau_<$) be a sign condition in
${\rm Sign}_{\ell-1}(\{X_1,\ldots,X_m\})$ in which some equation $X_i=0$ is replaced
by $X_i>0$ (respectively, by $X_i<0$).
Since functions in ${\mathcal P}_{\tau}$ are defined in $U_{\tau} \times V$ and
do not depend on $X_i$, they are also defined in the larger domain
$(U_{\tau} \cup U_{\tau_>} \cup U_{\tau_<}) \times V$.
It follows that for $S$ given by (\ref{eq:s}), the set $S'(\bar\eps)$ is well-defined
and is homotopy equivalent to $S$.

The rest of the proof is the same as the proof of Theorem~\ref{the:main2}.
\end{proof}

\section{On a conjecture of Benedetti and Risler}\label{sec:fewnomials}

In their book \cite{Risler}, Benedetti and Risler introduced an axiomatic definition
of the {\em complexity} of a polynomial, and consequently, of a semi-algebraic set.
They conjectured that there is a finite
number of homeomorphism types of semi-algebraic sets of a given complexity.
Benedetti and Risler  also formulated their conjecture for two particular cases:
the number of monomials and {\em additive complexity}.
In these cases the conjecture was proved independently by van den Dries \cite{Dries} and
Coste \cite{CosteFew} using o-minimality theory, without producing any explicit
upper bounds on the number of homeomorphism types as functions of complexity.

In this section we prove a weaker, but effective, versions of the two
mentioned particular cases of the conjecture.
Namely, we give explicit, single exponential
upper bounds on the number of possible {\em homotopy types} of semi-algebraic sets.
Our proofs are strongly influenced by arguments of van den Dries in
$\S$~3, Chapter~9 of \cite{Dries}.

\subsection{Homotopy types of sets defined by fewnomials}

We first prove a single exponential upper bound on the
number of homotopy types of
semi-algebraic subsets of $\Real^m$ defined by a family of polynomials
having in total at most $r$ monomials.

More precisely, let ${\mathcal M}_{m,r}$ be the family of ordered lists
${\mathcal P}=(P_1, \ldots , P_s)$ of polynomials $P_i \in \Real[X_1,\ldots,X_m]$,
with the total number of monomials in all polynomials in ${\mathcal P}$
not exceeding $r$.
Recall Definition~\ref{def:hom_lists}, of the homotopy type of an ordered
list of polynomials.

\begin{theorem}
\label{the:fewnomials}
The number of different homotopy types of ordered lists in ${\mathcal M}_{m,r}$
does not exceed
\begin{equation}\label{eq:fewnomials}
2^{O(mr)^4}.
\end{equation}
In particular, the number of different homotopy types of semi-algebraic sets defined
by a fixed formula $\phi_{\mathcal P}$, where ${\mathcal P}$ varies over ${\mathcal M}_{m,r}$,
does not exceed (\ref{eq:fewnomials}).
\end{theorem}

\begin{proof}
Observe that the function
$$|X|^Y= {\rm e}^{(Y \ln |X|)}$$
is Pfaffian of a constant order and degree in the domain
$\{ (x,y) \in \Real^2\ |\> x \neq 0 \}$.

We first prove the theorem for semi-algebraic sets defined in one fixed octant
$\Real_{\tau}^{m}$ for some sign condition $\tau \in {\rm Sign}(\{X_1, \ldots ,X_m \})$
(see Definition~\ref{def:octant}).

Let $\{i_1,\ldots,i_t\} \subset \{1,\ldots,n\}$ be such that,
$\tau(X_{i_j}) \in \{-1,1\}, 1 \leq j \leq t$, and
$\tau(X_{j}) =0 , j \in \{1,\ldots,n\}\setminus \{i_1,\ldots,i_t\}$.

Consider the family of $2^{r}$ Pfaffian functions
\begin{equation}\label{eq:pfaff}
\sum_{1 \le j \le r}\pm Z_j|X_{i_{1}}|^{Y_{ji_{1}}} \cdots |X_{i_{t}}|^{Y_{ji_{t}}}
\end{equation}
in $\Real^{m}_{\tau} \times \Real^{tr +r}$.
Observe that if
$$P:= \sum_{1 \le j \le r} q_jX_{i_{1}}^{\alpha_{ji_{1}}} \cdots X_{i_t}^{\alpha_{ji_{t}}}
\in \Real [X_{i_{1}}, \ldots ,X_{i_t}]$$
is a polynomial having $r$ monomials defined over
$\Real^{m}_{\tau}$, then there is a function
$$F_{\tau}(X_{i_{1}}, \ldots ,X_{i_t}, Y_{1i_{1}}, \ldots , Y_{ri_t}, Z_1, \ldots ,Z_r)$$
in the family (\ref{eq:pfaff}) such that
$$P=F_{\tau}(X_{i_{1}}, \ldots ,X_{i_t}, \alpha_{1i_{1}}, \ldots ,
\alpha_{ri_t}, q_1, \ldots ,q_r),$$
where the choice of the function is determined by parities
(values mod 2) of integers in the sequence
$\alpha_{ji_{1}}, \ldots ,\alpha_{ji_{t}}$.
For example, if $\tau(X_{i_c}) = -1$,
and $\alpha_{ji_{c}}$ is odd (respectively even),
then the term $|X_{i_{c}}|^{Y_{ji_{c}}}$
contributes the factor $(-1)$ (respectively 1) to the coefficient of the monomial
$Z_j|X_{i_{1}}|^{Y_{ji_{1}}} \cdots |X_{i_{t}}|^{Y_{ji_{t}}}$.

It follows that all polynomials having $r_0$ monomials defined in $\Real^{m}_{\tau}$
can be divided into $2^{r_0}$ disjoint classes so that each class is {\em represented} by a
Pfaffian function in the family (\ref{eq:pfaff}).

Consider an ordered list ${\mathcal P}_{\tau} = (P_1,\ldots,P_s)$ of polynomials
$P_i \in {\Real [X_{i_1},\ldots,X_{i_t}]}$.
Assume that the total number of monomials appearing
in all polynomials $P_1, \ldots ,P_s$ is $r$.
Let ${\mathcal F}_{\tau} = (F_{\tau 1},\ldots,F_{\tau s})$ be the list of Pfaffian functions
representing respective polynomials in ${\mathcal P}_{\tau}$.
We say that ${\mathcal F}_{\tau}$ {\em represents} ${\mathcal P}_{\tau}$.

Let
$$\pi_{\tau}:\> \Real^{m}_{\tau} \times \Real^{tr +r} \to \Real^{tr +r}$$
be the projection map.
According to Theorem \ref{the:main2}, for each list ${\mathcal F}_{\tau}$ there exists
a finite set $A_{{\mathcal F}_{\tau}} \subset \Real^{tr +r}$ with the cardinality
not exceeding
\begin{equation}\label{eq:bound}
2^{O(mr)^4},
\end{equation}
such that for every
$\y \in \Real^{tr +r}$ there exists $\z \in A_{{\mathcal F}_{\tau}}$ such that for
every ${\mathcal F}_{\tau}$-semi-Pfaffian set $S$, the set $\pi_{\tau}^{-1}(\y) \cap S$
is homotopy equivalent to $\pi_{\tau}^{-1}(\z) \cap S$.
In particular, (\ref{eq:bound}) is an upper bound on the number
of homotopy types of lists ${\mathcal P}_{\tau}$ which are represented by
${\mathcal F}_{\tau}$.
Since there are $2^r$ different lists ${\mathcal F}_{\tau}$, the number of homotopy
types of all lists ${\mathcal P}_{\tau}$
is at most (\ref{eq:bound}) multiplied by $2^r$.

Now we consider the general case of a semi-algebraic set defined in
$\Real^m$.
The proof will follow the same pattern as the just described simpler case of a set defined in
an octant $\Real^{m}_{\tau}$.
Observe that
$$\Real^m = \bigcup_{\tau \in {\rm Sign}(\{ X_1, \ldots , X_m \})} \Real^{m}_{\tau}.$$
Each ordered list ${\mathcal P}= (P_1,\ldots,P_s) \subset \Real[X_1, \ldots, X_n]$,
having in total $r$ monomials, and defined in $\Real^m$,
has its representative list ${\mathcal F}_{\tau} = (F_{\tau 1},\ldots,F_{\tau s})$ of
Pfaffian functions (among $2^r$ lists) when restricted on the octant
$\Real^{m}_{\tau}$.
Thus, each ${\mathcal P}$, defined in $\Real^m$, is {\em represented} by a family of lists
$\{{\mathcal F}_{\tau}|\> \tau \in {\rm Sign}(\{X_1, \ldots ,X_m \}) \}$ of Pfaffian
functions, corresponding to all $3^m$ octants.
For a given ${\mathcal P}$, the corresponding family of lists
$\{{\mathcal F}_{\tau}|\> \tau \in {\rm Sign}(\{X_1, \ldots ,X_m \}) \}$
is determined by the parities (values mod 2) of
integers in the total ordered sequence of powers of all $r$ monomials appearing in the
polynomials in ${\mathcal P}$.
Since the number of different sequences of parities of $r$ monomials
in $m$ variables is $2^{mr}$,
the total number of families we need to consider is bounded by $2^{mr}$.

For a Boolean formula $\phi$ with atoms
$\{a_i,b_i,c_i \mid 1 \leq i \leq s\}$ denote by $\phi_{{\mathcal F}_{\tau}}$
the formula obtained from $\phi$ by replacing
for each $i,\ 1\leq i \leq s $, the atom $a_i$
(respectively, $b_i$ and $c_i$) by
$F_{\tau i}= 0$ (respectively, $F_{\tau i} > 0$ and $F_{\tau i} < 0$).
Let $S(\phi_{{\mathcal F}_{\tau}})$ be the semi-Pfaffian set defined by
$\phi_{{\mathcal F}_{\tau}}$.

Let
$$\pi:\> \Real^{m} \times \Real^{mr +r} \to \Real^{mr +r}$$
be the projection map.
According to Theorem~\ref{the:main3}, for each fixed family of lists
$\{{\mathcal F}_{\tau}|\> \tau \in {\rm Sign}(\{X_1, \ldots ,X_m \}) \}$,
there exists a finite set $A \subset \Real^{mr +r}$
with cardinality not exceeding (\ref{eq:bound}), such that for every
$\y \in \Real^{mr +r}$ there exists $\z \in A$ such that for
every Boolean formula $\phi$ and
$$S:= \bigcup_{\tau \in {\rm Sign}(\{ X_1, \ldots , X_m \})}
S(\phi_{{\mathcal F}_{\tau}}),$$
the set $\pi_{\tau}^{-1}(\y) \cap S$ is homotopy equivalent to $\pi_{\tau}^{-1}(\z) \cap S$.
Hence (\ref{eq:bound}) is also an upper bound on the number of homotopy types of
lists ${\mathcal P}$ represented by the fixed family
$\{{\mathcal F}_{\tau}|\> \tau \in {\rm Sign}(\{X_1, \ldots ,X_m \}) \}$
of Pfaffian functions.
It follows that the number of all homotopy types of lists ${\mathcal P} \in
{\mathcal M}_{m,r}$ is at most (\ref{eq:bound}) multiplied by $2^{mr}$
since there are $2^{mr}$ different families of lists to consider.
\end{proof}

\subsection{Homotopy types of sets with bounded additive complexity}

\begin{definition}\label{def:additive}
A polynomial $P \in \Real[X_1, \ldots ,X_m]$ has {\em additive complexity
at most $a$} if there are polynomials $Q_1, \ldots , Q_a \in \Real[X_1, \ldots ,X_m]$
such that
\begin{itemize}
\item[(i)]
$Q_1=a_1X_{1}^{\alpha_{11}} \cdots X_{m}^{\alpha_{1m}} +
b_1X_{1}^{\beta_{11}} \cdots X_{m}^{\beta_{1m}}$,\\
where $a_1, b_1 \in \Real$, and
$\alpha_{11}, \ldots ,\alpha_{1m}, \beta_{11}, \ldots , \beta_{1m} \in \N$;

\item[(ii)]
$Q_j=a_jX_{1}^{\alpha_{j1}} \cdots X_{m}^{\alpha_{jm}}
\prod_{1 \le i \le j-1}Q_{i}^{\gamma_{ji}} +
b_jX_{1}^{\beta_{j1}} \cdots X_{m}^{\beta_{jm}}\prod_{1 \le i \le j-1}Q_{i}^{\delta_{ji}}$,\\
where $1 < j \le a$, $a_j, b_1j \in \Real$, and
$\alpha_{j1}, \ldots ,\alpha_{jm}, \beta_{j1}, \ldots , \beta_{jm}
\gamma_{ji}, \delta_{ji} \in \N$ for $1 \le i <j$;

\item[(iii)]
$P= cX_{1}^{\zeta_{1}} \cdots X_{m}^{\zeta_{m}}\prod_{1 \le j \le a}Q_{j}^{\eta_{j}}$,\\
where $c \in \Real$, and $\zeta_1, \ldots , \zeta_m, \eta_1, \ldots ,\eta_a \in \N$.
\end{itemize}
\end{definition}

In other words, $P$ has additive complexity at most $a$ if, starting with variables
$X_1, \ldots ,X_m$ and constants in $\Real$, and applying additions and multiplications,
a formula representing $P$ can be obtained using at most $a$ additions
(and an unlimited number of multiplications).

\begin{example}
The polynomial $P:=(X+1)^d \in \Real[X]$ with $0<d \in \Z$,
has $d+1$ monomials when expanded but additive complexity at most 1.
\end{example}

Let ${\mathcal A}_{m,a}$ be the family of ordered lists
${\mathcal P}=(P_1, \ldots , P_s)$ of polynomials $P_i \in \Real[X_1,\ldots,X_m]$,
with the additive complexity of every $P_k$ not exceeding $a_k$, and
$a=\sum_{1 \le k \le s}a_k$.

\begin{theorem}\label{the:additive}
The number of different homotopy types of ordered lists in ${\mathcal A}_{m,a}$
does not exceed
\begin{equation}\label{eq:additive}
2^{O((m+a)a)^4}.
\end{equation}
In particular, the number of different homotopy types of semi-algebraic sets defined
by a fixed formula $\phi_{\mathcal P}$, where ${\mathcal P}$ varies over ${\mathcal A}_{m,a}$,
does not exceed (\ref{eq:additive}).
\end{theorem}

\begin{proof}
Fix an ordered list ${\mathcal P} \in {\mathcal A}_{m,a}$.
For each polynomial $P_k \in {\mathcal P}$, $1 \le k \le s$,
consider the sequence of polynomials
$Q_{k1}, \ldots, Q_{ka_k}$ as in Definition~\ref{def:additive}, so that
$$
P_k:=c_kX_{1}^{\zeta_{k1}} \cdots X_{m}^{\zeta_{km}}\prod_{1 \le j \le a_k}Q_{kj}^{\eta_{kj}}.
$$
Introduce $a_k$ new variables
$Y_{k1}, \ldots ,Y_{ka_k}$.
Fix a semi-algebraic set $S \subset \Real^m$,
defined by a formula $\phi_{\mathcal P}$.
Consider
the semi-algebraic set $\widehat S$, defined by the conjunction of $a$ 3-nomial equations
obtained from equalities in (i), (ii) of Definition~\ref{def:additive} by replacing
$Q_j$ by $Y_{kj}$
for all $1 \le k \le s$, $1 \le j \le a_k$, and the formula $\phi_{\mathcal P}$
in which every occurrence of an atomic formula
of the kind $P_k \ast 0$, where $\ast \in \{ =, >, < \}$, is replaced by the formula
$$
c_kX_{1}^{\zeta_{k1}} \cdots X_{m}^{\zeta_{km}}\prod_{1 \le j \le a_k}Y_{kj}^{\eta_{kj}}
\ast 0.
$$
Note that $\widehat S$ is a
semi-algebraic set in $\Real^{m+a}$ defined by a formula
$\widehat \phi_{\widehat {\mathcal P}}$, where
$\widehat {\mathcal P} \in {\mathcal M}_{m+a, s+3a}$.
According to Theorem~\ref{the:fewnomials}, the number of different homotopy types of
ordered lists in ${\mathcal M}_{m+a, s+3a}$ does not exceed (\ref{eq:additive}).

Let $\rho:\> \Real^{m+a} \to \Real^m$ be the projection map on the subspace of coordinates
$X_1, \ldots ,X_m$.
It is clear that the restriction $\rho_{\widehat S}:\> \widehat S \to S$ is a homeomorphism.
Therefore,
if two lists ${\mathcal P},\ {\mathcal Q}\ \in {\mathcal A}_{m,a}$ are {\em not}
homotopy equivalent, then the corresponding lists
$\widehat {\mathcal P},\ \widehat {\mathcal Q}\ \in {\mathcal M}_{m+a, s+3a}$
are also not homotopy equivalent.
It follows, that the number of different homotopy types of lists in
${\mathcal A}_{m,a}$ does not exceed the number of different homotopy types of
lists in ${\mathcal M}_{m+a, s+3a}$, and therefore does not exceed (\ref{eq:additive}).
\end{proof}

\begin{remark}
A polynomial $P \in \Real[X_1, \ldots ,X_m]$ is said to have {\em rational additive complexity
at most $a$} if there are {\em rational functions}
$Q_1, \ldots , Q_a \in \Real(X_1, \ldots ,X_m)$ satisfying conditions (i), (ii), and (iii)
with $\N$ replaced by $\Z$.
For example, the polynomial $X^d+ \cdots + X+1 =(X^{d+1}-1)/(X-1)\in \Real[X]$
with $0<d \in \Z$, has rational additive complexity at most 2.
Define the family of ordered lists of polynomials ${\mathcal A}_{m, s, a}$ as above
but interpreting $a$
as the sum of rational additive complexities.
According to \cite{Dries} there is a finite the number of different homeomorphism types
of semi-algebraic sets defined by $s$ polynomials in $\Real [X_1, \ldots , X_m]$, with the
sum of additive complexities at most $a$.
We conjecture that the number of different homotopy types of lists in
${\mathcal A}_{m, s,a}$ does not exceed
$$
2^{(ma)^{O(1)}}.
$$
\end{remark}

\section{Metric upper bounds}
\label{sec:metric}
In this section we consider semi-algebraic sets defined by
polynomials with integer coefficients, and use the technique from previous sections to
obtain some ``metric'' upper bounds related to homotopy types.

Let $V \subset {\R}^m$ be a ${\mathcal P}$-semi-algebraic set,
where
${\mathcal P} \subset \Z [X_1, \ldots , X_m]$.
Let
for each $P \in {\mathcal P}$,
$\deg (P) <d$, and
the maximum of the absolute values of coefficients in $P$ be less than some
constant $M$,  $0 < M \in \Z$.
For $a > 0$ we denote by $B_m(0,a)$
the open ball of radius  $a$ in $\R^m$ centered at the origin.

The following proposition is well-known
(see, for instance, Theorem 4.1.1 in \cite{BPR95}).
\begin{proposition}
[\cite{V,GrV,HRS,BPR95}]\label{prop:ball}
There exists a constant $c > 0$, such that for any
$R > M^{d^{cm}}$, and
for any connected component $W$ of $V$ the intersection
$W \cap B_m(0,R)) \neq \emptyset$,
and $W \subset B_m(0,R)$ if $W$ is bounded.
\end{proposition}
Notice, that there
is no dependence in the bound on the cardinality of the family
${\mathcal P}$.
On the other hand, the dependence on the absolute values
of integer coefficients is essential,
as shown by an example of the semi-algebraic set defined by the system of equations
$$X=Y,\> X=(1- \eps)Y +1,$$
in which all coefficients are $\Theta (1)$,
independently of any small real $\eps >0$.

We prove the following generalization of Proposition~\ref{prop:ball}.

\begin{theorem}\label{the:ball}
There exists a constant $c > 0$, such that for any
$R_1 > R_2 > M^{d^{cm}}$ we have,
\begin{enumerate}
\item
$V \cap B_m(0,R_1) \simeq V \cap B_m(0,R_2)$, and
\item
$V \setminus B_m(0,R_1) \simeq V \setminus B_m(0,R_2)$.
\end{enumerate}
\end{theorem}

In the proof of Theorem \ref{the:ball} it will be convenient to use
{\em infinitesimals} instead of sufficiently small (or large) elements of
the ground real closed field $\R$.
We do this by considering non-archimedean extensions of $\R$.
More precisely, denote by
$\R\langle \varepsilon\rangle$  the real closed field of algebraic
Puiseux series in $\varepsilon$ with coefficients in $\R$
(see, for instance, \cite{BPRbook} for more details about Puiseux series).
The sign of a Puiseux series in $\R\langle \varepsilon\rangle$
agrees with the sign of the coefficient of the lowest degree term in $\varepsilon$.
This induces a unique order on $\R\langle \varepsilon\rangle$ which
makes $\varepsilon$
infinitesimal: $\varepsilon$ is positive and smaller than
any positive element of $\R$.
When $a \in \R\langle \varepsilon \rangle$ is bounded from above by an element of $\R$,
the symbol $\lim_\varepsilon(a)$ denotes the constant term of $a$, obtained by
substituting 0 for $\varepsilon$ in $a$, clearly $\lim_\varepsilon(a) \in \R$.
We will also denote by $\R\langle\bar\varepsilon\rangle$ the field
$\R\langle\varepsilon_1\rangle \cdots \langle \varepsilon_s\rangle$,
where $\varepsilon_{i+1}$ is infinitesimal with respect to
$\R\langle\varepsilon_1\rangle \cdots \langle \varepsilon_i\rangle$
for every $0 \le i < s$.

\begin{proof}[of Theorem \ref{the:ball}]
We will only prove the homotopy equivalence (1), the proof of (2)
being similar.

Let $\phi$ be the ${\mathcal P}$-formula defining $V$.
Let $S \subset \R^{m+1}$ be the set defined by the formula,
$\phi \wedge (X_1^2 +\cdots + X_m^2 - Y^2 \leq 0)$, and let
$\pi: \R^{m+1} \rightarrow \R$ be the projection on the $Y$ coordinate.
We will follow the notations introduced in the
proof of Theorem \ref{the:main}, but in the definition of
$S' = S'(\bar\eps)$ (cf. Definition \ref{def:S'})
we let $1 \gg \frac{1}{\omega} \gg \bar\eps > 0$ be infinitesimals and let $\R'$
be the field of Puiseux series, $\R\langle \frac{1}{\omega}\rangle \langle\bar\eps\rangle$.
The set $S'$ is then a semi-algebraic subset of $\R'^{m+1}$.

Consider now the set $G(S,\bar\eps) \subset \R'$
(cf. Definition \ref{def:criticalvalues}).
It follows as a consequence of the complexity analysis of
efficient quantifier elimination algorithms
(see \cite{GrV,BPRbook}),  that
$G(S,\bar\eps)$ is a finite subset of $\R'$
consisting of roots in $\R'$ of the univariate polynomials
$h \in \Z[\omega,\bar\eps][Y]$, whose degrees do not exceed $d^{O(k)}$, and
whose coefficients are integers with absolute values not exceeding $M^{d^{O(k)}}.$
Now suppose that $\alpha \in \R'$, $h(\alpha) =  0$ and $\alpha$ is bounded
from above by an element of $\R$. Then, $\lim_{\frac{1}{\omega}}(\alpha) \in \R$ is
a root of some polynomial in $\Z[Y]$ of degree not exceeding
$d^{O(k)}$, with absolute values of coefficients at most $M^{d^{O(k)}}$.
It follows that $|\alpha| \leq M^{d^{O(k)}}$.

Now let $G_b(S,\bar\eps) \subset G(S,\bar\eps)$
be the set of all elements of $G$ which are bounded from above by some element of $\R$.
Let $R = \max_{\alpha \in G_b} | \alpha |$,
then $R \leq  M^{d^{O(k)}}$.
Moreover, if elements $R_1,R_2 \in \R$, are both greater than $R$,
and belong to the same connected component of $\R' \setminus G(S,\bar\eps)$, then
$\pi_{S'}^{-1}(R_1) \simeq \pi_{S'}^{-1}(R_2).$
Since $\pi_S^{-1}(y) \simeq V$ for all sufficiently large $y \in \R$,
and $\pi_{S'}^{-1}(y) \simeq \pi_S^{-1}(y)$ for all $y \in \R$, it follows that
\[
\pi_{S'}^{-1}(R_1) \simeq \pi_{S'}^{-1}(R_2)\simeq \pi_S^{-1}(R_1) \simeq
\pi_S^{-1}(R_2) \simeq V.
\]
This proves the theorem.
\end{proof}

It is well-known that a semi-algebraic set has a {\em local conic structure}
(see, e.g., \cite{CosteSem}).
In particular, a semi-algebraic set is locally contractible.
The following theorem gives a quantitative version of the latter statement.

\begin{theorem}
Let $V \subset {\R}^m$ be a ${\mathcal P}$-semi-algebraic set,
with ${\mathcal P} \subset \Z [X_1, \ldots , X_m]$ and $0 \in V$.
Let $\deg (P) <d$ for each $P \in {\mathcal P}$, and
the maximum of absolute values of coefficients of $P \in {\mathcal P}$
be less than $M$, $0 < M \in \Z$.
Then, there exists a constant $c >0$ such that
for every $0 < r < M^{-d^{cm}}$ the set $V \cap B_m(0,r)$ is contractible.
\end{theorem}

\begin{proof}
Similar to the proof of Theorem~\ref{the:ball}.
\end{proof}

\affiliationone{
   Saugata Basu\\
   School of Mathematics,
    Georgia Institute of Technology, Atlanta, GA 30332\\
   USA
   \email{saugata.basu@math.gatech.edu}}
\affiliationtwo{
   Nicolai Vorobjov\\
   Department of Computer Science, University of Bath, Bath
   BA2 7AY, England\\
   UK
   \email{nnv@cs.bath.ac.uk}}

\begin{abstract}
In this paper we prove a single exponential upper bound on the number
of possible homotopy types of the fibres of a Pfaffian map, in terms
of the format of its graph.
In particular
we show that if a semi-algebraic set $S \subset {\R}^{m+n}$, where $\R$ is
a real closed field, is defined by a
Boolean formula with $s$
polynomials of degrees less than $d$, and
$\pi:\ {\R}^{m+n} \to {\R}^n$ is the projection
on a subspace, then the number of different homotopy types of fibres of
$\pi$ does not exceed
$s^{2(m+1)n}(2^m nd)^{O(nm)}$.
As applications of our main results
we prove single exponential bounds on the
number of homotopy types of semi-algebraic sets defined by fewnomials,
and by polynomials with bounded additive complexity.
We also prove single exponential
upper bounds on the radii of balls guaranteeing local contractibility
for semi-algebraic sets defined by polynomials with integer coefficients.
\end{abstract}

\end{document}